\documentclass[11pt]{article}
\usepackage[utf8]{inputenc}
\usepackage[english]{babel}
\usepackage{amsfonts}
\usepackage{amssymb,enumitem}
\usepackage{amsmath,amsthm}
\usepackage{booktabs}
\usepackage{graphicx}
\usepackage{titlesec}
\usepackage{adjustbox}

\usepackage{multirow}
\usepackage{caption}
\usepackage{cite}
\theoremstyle{definition}
\newtheorem{definition}{Definition}
\newtheorem{remark}{Remark}[section]

\newtheorem{theorem}{Theorem}
\newtheorem{lemma}{Lemma}
\newtheorem{prop}{Proposition}

\newcommand{\R}{\mathbb R}
\newcommand{\N}{\mathbb N}

\newcommand{\x}[2]{x_{#1}^{#2}}
\newcommand{\e}[2]{e_{#1}^{#2}}
\newcommand{\hx}[2]{\hat x_{#1}^{#2}}

\newcommand{\G}{\mathcal G}
\newcommand{\Z}{\mathcal Z}

\newcommand{\K}{\mathcal K}
\newcommand{\E}{\mathcal E}
\newcommand{\V}{\mathcal V}

\newcommand{\Rn}{\mathbb R^{n}}
 
\newcommand{\Rnn}{\mathbb R^{n\times n}}

\DeclareMathOperator{\diag}{diag}

\DeclareMathOperator{\Cov}{Cov}

\begin {document}
\title{Distributed Fixed Point Method for Solving Systems of Linear Algebraic Equations}
\author{Du\v{s}an Jakoveti\'c\footnote{Department of Mathematics and Informatics, Faculty of Sciences, University of Novi Sad, Trg Dositeja
Obradovi\'ca 4, 21000 Novi Sad, Serbia, Email: dusan.jakovetic@dmi.uns.ac.rs, natasak@uns.ac.rs, natasa.krklec@dmi.uns.ac.rs, greta.malaspina@dmi.uns.ac.rs.}, Nata\v{s}a Kreji\'c\footnotemark[1], Nata\v{s}a Krklec Jerinki\'c\footnotemark[1], \\ Greta Malaspina\footnotemark[1], Alessandra Micheletti\footnote{Department of Environmental Science and Policy - ESP, Università degli Studi di Milano, Via Saldini 50, 20133 Milano, Italy, Email: alessandra.micheletti@unimi.it}}
\date{}
\maketitle

\begin{abstract}
We present a class of iterative fully distributed fixed point methods to solve a system of linear equations, such that each agent in the network holds one of the equations of the system. Under a generic directed, strongly connected network, we prove a convergence result analogous to the one for fixed point methods in the classical, centralized, framework: the proposed method converges to the solution of the system of linear equations at a linear rate. We further explicitly quantify the rate in terms of the linear system and the network para\-meters. Next, we show that the algorithm provably works under time-varying directed networks provided that the underlying graph is connected over bounded iteration intervals, and we establish a linear convergence rate for this setting as well. A set of numerical results is presented, demonstrating practical benefits of the method over existing alternatives.\\

\end{abstract}

\section{Introduction}\label{derivation}
The problem we consider is \begin{equation}\label{eq:ls}
Ay=b
\end{equation}
where $A=[a_{ij}] \in\Rnn$ and $b=[b_i] \in\Rn$ are given, and $y\in\Rn$ is the vector of the unknowns. The matrix $A$ is assumed to be nonsingular, so that the problem has a unique solution. We also assume that the problem needs to be solved in a distributed computational framework determined by a set of connected computational nodes which can communicate through a generic sequence of graphs. 
Let  $A_i\in\R^{1\times n}$ and $b_i\in\R$ be the $i$-th row of $A$ and the $i$-th component of $b$ respectively. It is assumed that each computational node $ i $ knows the corresponding $ A_i $ and $ b_i $ and that each node needs to obtain the solution $ y^* $ through an iterative, distributed algorithm. 

The considered problem is important as linear systems appear naturally in a number of applications. One important example of application is to Ordinary Kriging \cite{kriging,krige,matheron}, an optimal linear prediction technique of the expected value of a spatial random field $\Z(s),\ s\in\mathbb{R}^n$. Ordinary Kriging can be applied when the random field under study is isotropic and intrinsically stationary, that is, the expected value $E(\Z(s))=m$ is constant and the variance $Var(\Z(s)-\Z(s+h))=2\gamma (h)$ depends only on $h$. In this case the model parameter estimation relies on the solution of a linear system like (\ref{eq:ls}) (see equations (3.2.13)-(3.2.15) in \cite{kriging} and the example in Section 5). When the semivariogram $\gamma (h)$ of the random field has a sill, it can be assumed that there is a range $\bar{h}$ over which the covariance $Cov(\Z(s),\Z(s+h))=0$, when $|h|>\bar{h}$. In this case the matrix $A$ of the Ordinary Kriging linear system becomes sparse, since its elements are the estimates of $\gamma(h)$ and each sampled node of the random field $\Z$ needs to memorize only the information brought by its neighbours at a distance lower than $\bar{h}$ to estimate the model parameters. When the mean $m$ of the random field is known the problem simplifies into what is called \emph{Simple Kriging}. In Section 5 we will use Simple Kriging as an example of application of our method.

There is a vast literature devoted to solving systems of linear equations in the conventional centralized computational environment \cite{greenbaum, saad}, as well as a number of results that cover parallelization of classical iterative methods which are applicable to the case of fully connected distributed computational environment, \cite{fs1}. Our interest in this paper is the class of fixed point methods \cite{greenbaum, saad} and their extensions to the distributed framework, as de\-scribed above.
In other words, we develop a class of novel, fully distributed, iterative fixed point methods to solve \eqref{eq:ls}, wherein each node can exchange messages only with the ones in its neighborhood in the communication graph, and each node obtains the estimate of the solution $y^*$ of problem \eqref{eq:ls}.
It is well known that \eqref{eq:ls} can be transformed into an equivalent fixed point problem
\begin{equation}
y = My + d, \label{fp}
\end{equation}
and one can apply the Banach contraction principle and define the fixed point iterative method of the form $ y^{k+1} = M y^k + d$, for suitable choices of $M\in\Rnn$ and $d\in\Rn$ (see Section 2 for the details). \\
The sufficient and necessary condition for the convergence of such iterative sequence  is $ \rho(M) < 1, $ where $ \rho(M) $ is the spectral radius of $ M. $ Furthermore, a sufficient condition for the convergence of $ \{y^k\} $ is given by $ \|M\| < 1 $ for an arbitrary matrix norm $ \|\cdot\|. $  Clearly, there is a number of suitable ways to define the iterative matrix $ M $ in such way that either $ \rho(M) < 1 $ or $ \|M\| < 1 $ for many matrix classes,  like symmetric positive definite matrices, M-matrices, H-matrices, etc \cite{Berman}. Typical examples of this type of methods are the  Jacobi and Gauss - Seidel method as well as their modifications like Jacobi Overrelaxation (JOR), Successive Overelaxation (SOR), Symmetric Successive Overelaxation (SSOR) method and so on \cite{greenbaum, saad}. The convergence of fixed point methods is linear and the convergence factor is determined by the spectral radius or the norm of $ M. $ The main idea of relaxation methods is to introduce a parameter that reduces the norm (or the spectral radius) of the corresponding iterative matrix and ensures faster convergence. 

There is a rich literature on parallelization of fixed point iterative meth\-ods, where the computational nodes communicate in an all-to-all fashion \cite{fs1}, \cite{fs2}, \cite{fm}, \cite{bertsekas}. In the case of very large dimension one needs to split the com\-pu\-ta\-tion\-al effort between different nodes to speed up the algorithm. In this type of computational environment, the total cost of solving the problem of interest is mainly dictated by the corresponding computational cost and the communication cost of exchanging messages between the parallelized nodes (processes) along iterations. Usually, major bottlenecks include waiting for the slowest node to complete an iteration, or latency incurred by the time to communicate a message. For this reason asynchronous methods, which allow for latency in communication and non\-uni\-form distribution of computational work, are also considered, \cite{fs2}. The methods of this type are convergent under different communication latency conditions \cite{fs2}.

The framework we consider in this paper for solving systems of linear equations, in more detail, assumes a network of computational nodes which communicate through a generic directed graph, which can depend on time. Thus the results in \cite{fs1, fm, fs2, bertsekas} are not applicable. The same framework is also considered in 
\cite{dls1,dls2,dls3,dls4,nedic}, and a survey of the methods is presented in \cite{survey}. The focus of these methods is to ensure convergence of the local approximations to the global solution, in the presence of time-varying com\-munica\-tion graphs. In the context of these algorithms, convergence is defined in two possible ways. In \cite{dls1, dls4} each node holds a local approximation of a subset of the variables and convergence of these local variables to the corresponding part of the solution is required. In \cite{dls2,dls3,nedic} every node contains a vector of the same size as the unknown vector of the linear system, and the convergence of each local vector to the full solution in ensured. We are interested in the second scenario.
The method presented in \cite{fs1} is applicable to a general problem of the type (\ref{eq:ls}) with loose restrictions on the matrix $ A $ and can be used to solve the linear least squares problem as well.

In this paper, we propose a novel distributed method to solve \eqref{eq:ls}, which we refer to as DFIX (Distributed Fixed Point). DFIX assumes the same computational framework as \cite{dls2,nedic,dls3} but differs significantly from the above mentioned methods. Namely, DFIX is derived starting from an as\-so\-ci\-at\-ed (centralized) fixed point method, rather than basing the derivation directly on the initial linear system. We extend the convergence theory of centralized fixed point methods to the distributed case in the sense of sufficient conditions. That is, we demonstrate that the sufficient condition $ \|M\|_{\infty} < 1 $ continues to work in the distributed environment. The main convergence result is completely analogous to the centralized case - given an iterative matrix with the infinity norm smaller than 1, the iterative sequence is convergent for an arbitrary starting point. The theory presented here thus covers a large class of linear systems. We prove linear convergence of DFIX under directed strongly connected networks and explicitly quantify the corresponding con\-ver\-gence factor in terms of network and linear system parameters. As detailed below, numerical simulations demonstrate ad\-van\-tages of DFIX over some state of the art methods.

With respect to the underlying graph, representing the connection among the computational agents, we consider both the case when the graph is fixed (i.e., the connectivity among the nodes is the same at any time during the execution of the algorithm) and the case when the network changes at every iteration. In the fixed graph case we prove that convergence holds if the network is strongly connected, while in the time-varying graph case we give suitable assumptions over the sequence of networks.
We prove that the time-independent case is a particular case of the time-varying case, but for the sake of clarity we first present and analyse the algorithm assuming the network is fixed, and then we generalize the analysis to the time-varying case.

Any system of linear equation \eqref{eq:ls} with symmetric matrix $A$ can be considered as the first order optimality condition of an unconstrained op\-ti\-miza\-tion problem with cost function $\frac{1}{2}x^tAx-b^tx.$ It is therefore of interest to compare the approach of solving \eqref{eq:ls} applying some distributed optimization method \cite{harnessing, diging, extra} to the minimization of the quadratic function $\frac{1}{2}x^tAx-b^tx$ with DFIX. We thus compare computational and communication costs of DFIX with the state of the art optimization method from \cite{harnessing} and show that the computational costs with DFIX are significantly lower,  while the com\-mu\-ni\-ca\-tion costs are comparable or go in favor of DFIX, depending on the connectivity of the underlying graph. Thus the numerical efficiency of DFIX is also shown. A comparison with the method from \cite{nedic} is also presented in Section 5, demonstrating the clear advantage of DFIX. 

This paper is organized as follows. Section 2 contains the description of the computational framework together with a brief overview of fixed point iterative methods that will be used further on. The method DFIX is defined and analysed in Section 3 for the fixed graph case. In Section 4 we present the time-varying case. Numerical results that illustrate theoretical analysis as well as an application of DFIX to a kriging problem are presented in Section 5. Some conclusions are drawn in Section 6. 

\section{Preliminaries}

  Let us first briefly recall the theory of fixed point iterative methods for systems of linear equations.  Given a generic \footnote{the relation between the method in \eqref{fpiter} and \eqref{eq:ls} is described further ahead} method of type (\ref{fp})
  \begin{equation}
x^{k+1} = M x^k + d,  \label{fpiter}
\end{equation}
   we know that the method is convergent if $ \rho(M) < 1, $ where we recall that $ \rho(M) $ is the spectral radius of $ M$, i.e., the largest eigenvalue of $M$ in modulus. This condition is both necessary and sufficient for convergence. Given any matrix norm $ \|\cdot\| $ one can also state the sufficient convergence condition as $ \|M\| < 1. $ There are many ways of transforming (\ref{eq:ls}) to the fixed point form (\ref{fp}), depending on the properties of $ A, $ with Jacobi and Gauss - Seidel methods, as well as their relaxation versions being the most studied methods. To fix the idea before defining the distributed method we recall here the Jacobi and Jacobi Overrelaxation, JOR, method, keeping in mind that we will consider a generic $ M $ in the next section. 
  
  Assume that $ A $ is a nonsigular matrix with nonzero diagonal entries. Using the splitting 
  $ A = D -  P,  $ with $ D $ being the diagonal matrix, $ D = \diag(a_{11},\ldots,a_{nn}), $  the Jacobi iterative method is defined by (\ref{fpiter}) with  
  $$ M = D^{-1}P:=M_{J}. $$  In other words, given $ d=D^{-1}b $ and denoting by $x^k = (x^k_1,\dots,x^k_n)$ the
estimate of solution to \eqref{eq:ls} at iteration $k$, the new iteration is defined by
  $$ x_{i}^{k+1} = -\frac{1}{a_{ii} }\sum_{j=1,j\neq i}^n a_{ij} x_{j}^k + d_i, \; i=1,\ldots,n. $$
  The method is linearly convergent for many classes of matrices, for example strictly diagonally dominant matrices, symmetric positive definite matrices etc \cite{greenbaum, saad}, and the rate of convergence is determined by $ \rho(M_{J}). $ To speed up convergence and extend the class of matrices for which the method is convergent, one can introduce the relaxation parameter $ \omega \in \mathbb{R}$ and define
  $$    M = \omega D^{-1}P + (1-\omega) I. $$ In other words, the JOR iteration is given by 
  \begin{equation} \label{JOR}
  x_{i}^{k+1} = (1-\omega) x_{i}^k - \frac{w}{a_{ii} }(\sum_{j=1,j\neq i}^n a_{ij} x_{j}^k + b_i), \; i=1,\ldots,n. 
  \end{equation} 
  If $ A $ is a symmetric positive definite matrix, the JOR method converges for 
  $$ \omega \in (0, \frac{2}{\rho(M_J)}), $$ see \cite{greenbaum, saad}.  
  
  Assuming that each node can communicate directly with every other node, the method can be applied in parallel and asynchronous manner and the convergence follows from the results of \cite{fs2,bertsekas}. 
 
Let us now define precisely the computational environment we consider.   Assume that the network of nodes is a directed network $ {\cal G} = ({\cal V},{\cal E}), $ where $ {\cal V} $ is the set of nodes and $ {\cal E} $ is the set of all edges,  i.e., all pairs $(i,j) $ of nodes where node $i$ can send information to node $j$  through a communication link. 
\begin{definition}
The graph $\G=(\V,\E)$ is \emph{strongly connected} if for every couple of nodes $i,j$ there exists an oriented path from $i$ to $j$ in $\G$. That is, if there exist ${s_1,\dots, s_l}$ such that $(i,s_1), (s_1, s_2), \dots, (s_l,j)\in\E.$
\end{definition}

 \noindent {\bf Assumption A1.}
The network  $ {\cal G} = ({\cal V},{\cal E}) $ is directed, strongly connected, with self-loops at every node. 
 
 \begin{remark}
 The case of undirected network $\G$ can be seen as the particular case of directed graph where $\G$ is symmetric. That is, $(i,j)\in\E$ if and only if $(j,i)\in\E.$ In this case, the hypothesis that $\G$ is strongly connected is equivalent to $\G$ connected.
 \end{remark}

 Let us denote by $ O_i $ the in-neighborhood of node $i$, that is, the set of nodes that can send information to node $i$ directly. Since the graph has self loops at each node, then $i\in O_i$ for every $i$. We associate with $ {\cal G} $ an $n \times n$ matrix $ W$, such that the elements of $ W $ are all nonnegative and each row sums up to one. More precisely, we assume the following.

 \noindent{\bf Assumption A2.} The matrix $ W \in \mathbb{R}^{n \times n} $ is row stochastic
 with elements $ w_{ij} $ such that
 $$ w_{ij} > 0 \mbox{ if } j\in O_i, \;  w_{ij} = 0 \mbox{ if }  j \notin O_i$$
 
 Let us denote by $w_{min}$ a  constant such that 
  all nonzero elements of $ W $ satisfy $   w_{ij} \geq w_{\min}>0. $ Under the previously stated assumptions we know that such constant exists. Moreover, we have $w_{\min} \in (0,1)$. 
 Therefore, for all elements of $ W $ we have 
 \begin{equation} \label{wmin} 
 w_{ij} \neq 0 \Rightarrow w_{ij} \geq w_{\min}. 
  \end{equation}
  The diameter of a network is defined as the largest distance between two nodes in the graph. Let us denote with $\delta$ the diameter of $\G.$

\section{DFIX method}

We consider now a generic fixed point method for solving \eqref{eq:ls} by the fixed point iterative method (\ref{fpiter}), 
with $\ M=[m_{ij}] \in\Rnn,\ d=[d_i] \in\Rn$ defined in such a way that node $i$ contains the $i$-th row $M_i\in\R^{1\times n}$ and $d_i\in\R.$ Moreover, we assume that the fixed point $y^*$ of \eqref{fp} is a solution of \eqref{eq:ls}. The algorithm is designed in such way that each node  has its own estimate of the solution $ y^*. $ Thus at iteration $ k $ each node $ i $ has its own estimate $ x_i^k \in \mathbb{R}^n $ with components $ x_{ij}^k, \; j=1,\ldots,n. $ 
The DFIX method is presented in the algorithm below.\\

\noindent{\bf Algorithm DFIX}
\begin{itemize}
\item[Step 0] Initialization: Set $k = 0$. Each node chooses $ x_i^0 \in \mathbb{R}^n. $ 
\item[Step 1] Each node $i$ computes 
 \begin{equation}\label{eq:step1}\begin{aligned}
&\hx{ii}{k+1}=\sum_{j=1}^n m_{ij}\x{ij}{k}+d_i,\\
& \hx{ij}{k+1} =\hx{ij}{k},\ i \neq j.  \\
\end{aligned}
\end{equation}
\item[Step 2] Each node $i$ updates its solution estimate
 \begin{equation}\label{eq:step2}\x{i}{k+1}=\sum_{j=1}^nw_{ij}\hx{j}{k+1}\end{equation} and sets 
 $ k = k+1. $
\end{itemize}

Notice that at Step 1 each node $ i $ updates only the $i$-th component of its solution estimate and leaves all other components unchanged, while in Step 2 all nodes perfom a consensus step \cite{consensus, graphs, touri} using the set of vector  estimates $ \hx{j}{k+1}. $ 
Defining the global variable at iteration $k$ as
\begin{equation*}
X^k=\left( x_1^k;  \ldots;  x_n^k \right)\in\R^{n^2}, 
\end{equation*}  
Algorithm DFIX can be stated in a condensed form using $ X^k $ and  the following notation
\begin{equation*}
\widehat M_i=
\left(\begin{array}{ccccc}
1& &  & & \\
 &\ddots & & &\\
m_{i1} &\dots &m_{ii}  &\dots & m_{in}\\
 & &  &\ddots &\\
 & &  & &1\\
\end{array}\right)\in\Rnn, \phantom{spa}
\widehat d_i=\left(\begin{matrix}
 0 \\
 \vdots\\
 d_i\\
 \vdots\\
 0 \\
\end{matrix}\right)\in\Rn.
\end{equation*}
More precisely, matrix $\widehat M_i$ has the $i$-th row equal to $M$, the rest of diagonal elements are equal to 1 and the remaining elements are equal to 0. Vector $\widehat d_i$ has only one nonzero element in the $i$-th row which is equal to $d_i$.  
Now, Step 1 can be rewritten as 
\begin{equation*}
\hx{i}{k+1}=\widehat{M}_i\x{i}{k}+\hat{d}_i,
\end{equation*} and we can rewrite the Steps 1-2  in matrix form as
\begin{equation}\label{eq:globalit}
X^{k+1}=(W\otimes I)(\mathcal MX^k+\hat d)
\end{equation}
where $\mathcal M=\diag\left(\widehat M_1,\dots,\widehat M_n\right)\in\R^{n^2\times n^2},$ $\widehat d=\left(\hat d_1;\dots;\hat d_n\right)\in\R^{n^2}$ and $\otimes $ denotes the Kronecker product of matrices. 
We remark here that equation \eqref{eq:globalit} is only theoretical, in the sense that since each agent has access only to partial information, the global vector $X^k$, the matrix $\mathcal M$ and the vector $\widehat d$ are not computed at any node. We derived equation \eqref{eq:globalit} to get a compact representation of Algorithm 1 and to use it in the convergence analysis.

The following theorem shows that for every $i \in \{1,\ldots,n\} $ the local sequence $\{\x{i}{k}\}$ converges to the fixed point $y^*$ of \eqref{fp}. Denote 
\begin{equation*}
X^*=\left(y^*; \ldots; y^*\right) \in\R^{n^2}.  
\end{equation*}

\begin{theorem}
Let Assumptions A1 and A2 hold, $\|M\|_\infty=\mu < 1 $ and let $\{X^k\}$ be a sequence generated by \eqref{eq:globalit}. 
There exists a constant $\tau <1$ such that for every $k$ the global error $E^{k}=X^k-X^*$ satisfies 
 \begin{equation}
\|E^{k+1}\|_{\infty}\leq\tau \|E^{k-\delta+1}\|_{\infty},
\end{equation}
where $\delta$ denotes the diameter of the underlying computational graph $  {\cal G}.$
\end{theorem}

\begin{proof}
Since $W$ is assumed to be row stochastic there holds 
$(W\otimes I)X^*=X^*$. Moreover, using the fact that 
$\hat d=(I \otimes I-\mathcal M)X^*$, we obtain the following recursion
\begin{equation}\label{novo1} E^{k+1}=(W\otimes I) \mathcal M E^k.
\end{equation}
Notice that $\|(W\otimes I) \mathcal M\|_{\infty} \leq 1$, so we have 
\begin{equation}\label{novo2}\|E^{k+1}\|_{\infty}\leq \|E^k\|_{\infty}.
\end{equation}
Now, denoting by $e_i^k$ the $i$-th block of $E^k$ (the local error corresponding to node $i$) and by $e_{ij}^k$ its $j$-th component, from \eqref{novo1} we obtain the following 
\begin{equation}\label{eq:errorcomp}
\e{ij}{k+1}=w_{ij}M_j\e{j}{k}+\sum_{s\neq j}w_{is}\e{sj}{k}.
\end{equation}
We prove the thesis by proving that if the distance between $j$ and $i$ in the graph is equal to $l$, then for every $k$
\begin{equation}\label{eq:distl}
|\e{ij}{k+1}|\leq\tau'\|E^{k-l+1}\|_{\infty},\ \text{for a constant}\ \tau'<1.
\end{equation}
We proceed by induction over the distance $l$.
If $l=1$, that is, if there is an edge from $j$ to $i$, then $w_{ij}\geq w_{min}>0$. By \eqref{eq:errorcomp} we get
\begin{equation*}\label{eq:diagerrderiv}
\begin{aligned}
|\e{ij}{k+1}|&\leq w_{ij}|M_j\e{j}{k}|+\sum_{s\neq j}w_{is}|\e{sj}{k}|
\leq 
w_{ij}\mu\|E^{k}\|_{\infty}+\|E^{k}\|_{\infty}\sum_{s\neq j}w_{is}\leq\\
&\leq\big(1-w_{ij}(1-\mu)\big)\|E^{k}\|_{\infty}\leq\big(1-w_{\text{min}}(1-\mu)\big)\|E^{k}\|_{\infty}, 
\end{aligned}
\end{equation*}
and defining $\tau' =\big(1-w_{\text{min}}(1-\mu)\big)<1,$ we get
\begin{equation}\label{eq:tau1}
|\e{ij}{k+1}|\leq\tau'\|E^{k}\|_{\infty}.
\end{equation}
Assume now that \eqref{eq:distl} holds for distance equal to $l-1$, and let us prove it for $l$. Let $(j, s_{l-1}, s_{l-2}, \dots, s_1, i)$ be a path of length $l$ from $j$ to $i$. In particular we have that $w_{is_{1}}>0$ and thus 
\begin{equation}\label{eq:s1}
|\e{ij}{k+1}|\leq w_{is_1}|\e{s_1j}{k}|+\sum_{s\neq s_1}w_{is}|\e{sj}{k}|.
\end{equation}
For each of the terms $|\e{sj}{k}|$ in the sum, by \eqref{novo2}, we have
\begin{equation}\label{eq:1}
|\e{sj}{k}|\leq\|E^k\|_{\infty}\leq\|E^{k-l+1}\|_{\infty}.
\end{equation}
Let us now consider the term $|\e{s_1j}{k}|$. Since $(j, s_{l-1}, s_{l-2}, \dots, s_1, i)$ is a path of length $l$ from $j$ to $i$ and the distance between $j$ and $i$ is equal to $l$, we have that the distance between $j$ and $s_1$ is equal to $l-1$ and therefore, by inductive hypothesis
\begin{equation}\label{eq:2}
|\e{s_1j}{k}|\leq\tau'\|E^{k-(l-1)}\|_{\infty} = \tau'\|E^{k-l+1}\|_{\infty},\ \text{for}\ \tau'<1.
\end{equation}
Replacing \eqref{eq:1} and \eqref{eq:2} in \eqref{eq:s1}, we get 
\begin{equation}
\begin{aligned}
|\e{ij}{k+1}|&\leq w_{is_1}\tau'\|E^{k-l+1}\|_{\infty}+\sum_{s\neq s_1}w_{is}\|E^{k-l+1}\|_{\infty}=\\
& = \left(1-w_{s_1j}(1-\tau')\right)\|E^{k-l+1}\|_{\infty}\leq\\
&\leq \left(1-w_{min}(1-\tau')\right)\|E^{k-l+1}\|_{\infty}
\end{aligned}
\end{equation}
and defining $\tau:=\left(1-w_{min}(1-\tau')\right)<1$ we get \eqref{eq:distl}. 
Now the thesis follows directly from the fact that the distance between any two nodes is smaller or equal than the diameter $\delta$ of the graph.

\end{proof}

\section{Time-varying Network}
The method discussed in the previous sections is valid only if the graph representing the communication among the agents is the same at each iteration. If some failure of the communication link between two agents occurs during the execution of the algorithm, the underlying network changes, and Theorem 1 does not apply anymore. To deal with these possible changes we consider the case where the network is given, possibly different, at each iteration. We extend DFIX to this framework and we give assumptions on the sequence of graphs that yield a convergence result analogous to Theorem 1. In particular we show that, in order to achieve convergence, strong connectivity is not necessary at any time.

Assume that a sequence of directed graphs $\{\mathcal G_k\}_k$ is given, such that $\mathcal G_k$ represents the network of nodes at iteration $k$. That is, at iteration $k$, each node can communicate with its neighbours in $\mathcal G_k$. The DFIX algorithm described by equations \eqref{eq:step1} and \eqref{eq:step2} can be applied in this case if we replace \eqref{eq:step2} with
\begin{equation}\label{eq:step2_time}
\x{i}{k+1}=\sum_{j=1}^nw_{ij}^k\hx{j}{k+1}
\end{equation}
where $W^k$ is the consensus matrix associated with the graph $\mathcal G_k$, that is, $W^k$ satisfies Assumption A2 with $\G = \G_k$. With this modification, the equation describing the global iteration becomes 
\begin{equation}\label{eq:globalit_time}
X^{k+1}=(W^k\otimes I)(\mathcal MX^k+\hat d).
\end{equation}
We will prove a convergence result for a class of sequences of graphs. We first present and analyze the assumptions on such sequence.
\begin{definition}
Given $\mathcal G_1, \mathcal G_2$ graphs with $\mathcal G_i=(\mathcal V, \mathcal E_i)$, the \emph{composition} of $\mathcal G_1$ and $\mathcal G_2$ is defined as $\mathcal G_2\circ \mathcal G_1 = (\mathcal V, \mathcal E)$ where
\begin{equation}
\mathcal E := \{(j,i)\in \mathcal V^2\ |\ \exists\ s\in\V\ \text{such that}\ (j,s)\in \mathcal E_1, (s,i)\in \mathcal E_2\}.
\end{equation}
\end{definition} 
That is, there is an edge from $j$ to $i$ in $\mathcal G_2\circ \mathcal G_1$ if we can find a path from $j$ to $i$ such that the first edge of the path is in $\mathcal G_1$ and the second edge is in $\mathcal G_2$. This definition can be extended to finite sequences of graphs of arbitrary length. 
\begin{remark}\label{edgescomposition}
Let us consider a generic set of graphs $\G_1,\dots,\G_{m}$.
It is easy to see that if for every index $j$ the graph $\G_j$ has self-loops at every node then the set of edges of the composition $\G_1\circ\dots\circ\G_{m}$ contains the set of edges of $\G_j$ for every $j$. In particular, if there exists an index $\hat\jmath\in\{1,\dots,m\}$ such that $\G_{\hat\jmath}$ is fully connected, then  $\G_1\circ\dots\circ\G_{m}$ is also fully connected.
\end{remark} 
\begin{definition}Given an infinite sequence of networks $\{\mathcal G_k\}_k$ and a positive integer $\bar m$, we say that the sequence is \emph{jointly fully} (respectively, \emph{strongly}) \emph{connected for sequences of length $\bar m$} if for every index $k$, the composition $\mathcal G_k\circ \mathcal G_{k+1}\circ\dots\circ\mathcal G_{k+\bar m-1}$ is fully (respectively, strongly) connected. 
\end{definition}

\begin{definition}Given an infinite sequence of networks $\{\mathcal G_k\}_k$ and two integers $\tau_0$, $l$, we say that the sequence is \emph{repeteadly jointly strongly connected with constants $\tau_0$, $l$}, if for every index $k$, the composition $\mathcal G_{\tau_0+kl}\circ \mathcal G_{\tau_0+kl+1}\circ\dots\circ\mathcal G_{\tau_0+(k+1)l}$ is strongly connected. 
\end{definition}

\begin{definition} Given two vertices $i, \ j$ we say that there is a \emph{joint path} of length $l$ from $i$ to $j$ in $\G_k,\dots, \G_{k+\bar m -1}$ if there exist $s_1,\dots, s_{l-1}$ such that $(i,s_1)\in \mathcal E_{k+\bar m -1}, \ (s_1,s_2)\in \mathcal E_{k+\bar m-2}, \dots, (s_{l-1},j)\in \mathcal E_{k+\bar m-l}$, and we say that $i, j$ have \emph{joint distance} $l$ in $\G_k,\dots, \G_{k+\bar m -1}$ if the shortest joint path from $i$ to $j$ is of length $l$. \end{definition}

\noindent Our analysis is based on the following assumption.\\

\noindent {\bf Assumption A3.}
$\{\mathcal G_k\}$ is a sequence of directed graphs, with self-loops at every node, jointly fully connected for sequences of length $\bar m$, for some positive integer $\bar m$.\\

\noindent The algorithm presented in \cite{dls2} works for time-varying network in a similar framework. Formally, the hypothesis on $\{\G_k\}$ in \cite{dls2} is the following.\\

\noindent {\bf Assumption A3'.}
$\{\mathcal G_k\}$ is a sequence of directed graphs, with self-loops at every node, jointly strongly connected for sequences of length $\bar p$, for some positive integer $\bar p$.\\

\noindent We show now that Assumptions A3 and A3' are equivalent, in the sense specified by Proposition 1. In the following, given an integer $m$, we denote with $\G^m$ the composition of $m$ copies of $\G$.
 \begin{lemma}
 If $\G$ is a directed strongly connected graph with self-loops at every node and diameter $\delta$, then $\G^\delta$ is fully connected.
 \begin{proof}
 By definition of composition we have that $(i,j)$ is an edge in $\G^{\delta}$ if and only if 
 \begin{equation}\label{compositioncondition}
 \exists s_1,\dots, s_{\delta-1}\in\V\ \text{such that}\  
 (i,s_1), (s_1, s_2), \dots, (s_{\delta-1},j) \in \G.
 \end{equation}
 We want to prove that for every $i,j\in\V$ a sequence of nodes $s_h$ as in \eqref{compositioncondition} exists.\\
 Since $\G$ is fully connected with diameter $\delta$, there exists a path in $\G$ from $i$ to $j$ of length $l\leq\delta$. That is, there exist a set of nodes $v_1, \dots, v_{l-1}$ such that 
 $(i,v_1), (v_1, v_2), \dots, (v_{l-1},j)$ are edges in $\G$ and therefore a sequence satisfying \eqref{compositioncondition} is given by 
 \begin{equation*}
 s_h = \begin{cases}
 v_h & h = 1:l-1 \\ j & h = l:\delta.
 \end{cases}
 \end{equation*}
 \end{proof}
 \end{lemma}

\begin{prop}
Let $\{\G_k\}$ be a sequence of graphs where, for each $k$, $\G_k = (\V, \E_k)$ is a directed graph with self-loops at every node. The following are equivalent:
\begin{enumerate}
\item[(1)] there exist $\tau_0, l\in\mathbb N$ such that $\{\G_k\}$ is repeatedly jointly strongly connected with constants $\tau_0, l$
\item[(2)] there exists $\bar p\in\mathbb N$ such that $\{\G_k\}$ is strongly connected for sequences of length $\bar p$
\item[(3)] there exists $\bar m\in\mathbb N$ such that $\{\G_k\}$ is fully connected for sequences of length $\bar m$
\end{enumerate}
\end{prop}
\begin{proof}
It is easy to see that $(2)\Rightarrow(1)$ with $\tau_0 = 0$ and $l=\bar p$ and since full connectivity clearly implies strong connectivity, we have that $(3)\Rightarrow(2)$ with $\bar p = \bar m.$\\
We now prove that $(1)\Rightarrow(2)$ with $\bar p = 2l.$ That is, we prove that if (1) holds, then for every index $s$ the composition
$\G_s\circ\dots\circ\G_{s+2l-1}$ is strongly connected.
Given an index $s$, we denote with $\bar r$ the remainder of the division of $(s-\tau_0)$ by $l$, we define $\bar h:=l^{-1}(s-\tau_0+l-\bar r)$. By definition of $\bar r$ and $\bar h$ and applying (1) with $k = \bar h$ we have that the graph
\begin{equation*}\begin{aligned}
H:&= \G_{s+l-\bar r}\circ\dots\circ\G_{s+2l-\bar r-1}=\\  
&= \G_{\tau_0+\bar h l}\circ\dots\circ\G_{\tau_0+(\bar h+1) l-1} 
\end{aligned}
\end{equation*}
is strongly connected and thus 
\begin{equation*}\begin{aligned}
&\G_s\circ\dots\circ\G_{s+2l-2} =  \G_s\circ\dots\circ\G_{s+l-\bar r-1}\circ H\circ \G_{s+2l-\bar r}\circ\dots\circ\G_{s+2l-1}
\end{aligned}\end{equation*}
is strongly connected. Since $2l-\bar r\in{l+1,\dots,2l}$ we have the thesis.\\
Finally, we prove that $(2)\Rightarrow(3)$.
Since the size of $\V$ is finite, there exists a finite number of graphs with vertices $\V$. In particular, there exists a finite integer $L$ equal to the number of strongly connected graphs with vertices $\V.$ We denote with $H_1, \dots H_L$ such graphs, with $\delta_j$ the diameter of $H^j$ and with $\bar\delta:= \max{\delta_j}.$
Given any index $k$, we consider $(\bar\delta-1)L+1$ sequences of length $\bar p$ as follows:
\begin{equation*}\begin{aligned}
&S_1 = \G_{k}\circ\G_{k+1}\dots\circ\G_{k+\bar p-1}\\
&S_2 = \G_{k+\bar p}\circ\G_{k+\bar p+1}\dots\circ\G_{k+2\bar p-1}\\
&\vdots\\
&S_{(\bar\delta-1)L+1} = \G_{k+(\bar \delta-1)L\bar p}\circ\G_{k+(\bar \delta-1)L\bar p+1}\dots\circ\G_{k+(\bar \delta-1)L\bar p+\bar p-1}.
\end{aligned}
\end{equation*}
Statement (2) implies that, for every $j\in\{1,\dots,(\bar\delta-1)L+1\},\ S_j\in\{H_1, \dots H_L\}$ and thus there exists an index $\hat\imath\in\{1,\dots,L\}$ such that at least $\bar\delta$ elements of $\{S_1,\dots, S_{(\bar\delta-1)L+1}\}$ are equal to $H_{\hat\imath}$. Using the fact that, by Lemma 1, $H_{\hat\imath}^{\delta_{\hat\imath}}$ is fully connected and Remark \ref{edgescomposition}, we have
\begin{equation*}\begin{aligned}
\G_{k}\circ\G_{k+1}\circ\dots\circ\G_{k+(\bar \delta-1)L\bar p+\bar p-1}  = S_1\circ\dots\circ S_{(\bar\delta-1)L+1}
\end{aligned}
\end{equation*}
fully connected, and thus (3) holds with $\bar m = (\bar \delta-1)L\bar p+\bar p$.\\

\end{proof}

To conclude the considerations on the sequence of networks we remark that, since we are assuming that the linear system \eqref{eq:ls} has unique solution and that each node contains exactly one row of the coefficient matrix, the $D$-connectivity hypothesis introduced in  \cite{nedic} is equivalent to Assumption A3' and thus, by Proposition 1, to Assumption A3.

\begin{theorem}
Assume that a sequence of networks $\{\mathcal G_k\}_k$ is given, satisfying Assumption A3, and that for every index $k$ the corresponding consensus matrix $W^k$ satisfies Assumption A2. Let $\{X^k\}$ be a sequence generated by \eqref{eq:globalit_time} with 
$\|M\|_\infty=\mu < 1 $.
There exists a constant $\sigma <1$ such that for every $k\in\N$ the global error $E^{k}=X^k-X^*$ satisfies 
 \begin{equation}
\|E^{k+1}\|_{\infty}\leq\sigma \|E^{k-\bar m+1}\|_{\infty},
\end{equation}
where $\bar m$ is the constant given by Assumption A3.
\begin{proof} We follow the proof of Theorem 1.
For every index $k$, the matrix $W^k$ is row stochastic and $\|(W^k\otimes I) \mathcal M\|_{\infty} \leq 1$, so we get
\begin{equation}\label{novo1_time} E^{k+1}=(W^k\otimes I) \mathcal M E^k.
\end{equation}
and \begin{equation}\label{novo2_time}\|E^{k+1}\|_{\infty}\leq \|E^k\|_{\infty}.
\end{equation}
For every node $i$, $j$ and for every iteration index $k$, we have
\begin{equation}
\e{ij}{k+1}=w_{ij}^kM_j\e{j}{k}+\sum_{s\neq j}w_{is}^k\e{sj}{k}.
\end{equation}
We now prove that if the joint distance between $j$ and $i$ in $\G_{k-\bar m+1}$, $\G_{k-\bar m+2}$, $\dots, \G_k$ is equal to $l$, then for every $k$
\begin{equation}\label{eq:distltime}
|\e{ij}{k+1}|\leq\sigma'\|E^{k-l+1}\|_{\infty},\ \text{for}\ \sigma'<1.
\end{equation}
We proceed by induction over the joint distance $l$.
If $l=1$, that is, if $w_{ij}^k>0$, proceeding as in the derivation of \eqref{eq:diagerrderiv} we get
\begin{equation*}
\begin{aligned}
|\e{ij}{k+1}|&\leq\big(1-w_{ij}^k(1-\mu)\big)\|E^{k}\|_{\infty}\leq\big(1-w_{\text{min}}(1-\mu)\big)\|E^{k}\|_{\infty}=:\sigma\|E^{k}\|_{\infty}.
\end{aligned}
\end{equation*}
We assume now that \eqref{eq:distltime} holds for distance equal to $l-1$ and we prove it for $l$. Let $(j, s_{l-1}, s_{l-2}, \dots, s_1, i)$ be a joint path of length $l$ from $j$ to $i$ in 
$\G_{k-\bar m+1}, \G_{k-\bar m+2},\dots, \G_k$
 In particular we have that $w_{is_{1}}^k>0$ and thus 
\begin{equation}\label{eq:s1time}
|\e{ij}{k+1}|\leq w_{is_1}|\e{s_1j}{k}|+\sum_{s\neq s_1}w_{is}|\e{sj}{k}|.
\end{equation}
Using the fact that $(j, s_{l-1}, s_{l-2}, \dots, s_1)$ is a joint path of length $l-1$ from $j$ to $s_1$ in 
$\G_{k-\bar m+1}, \G_{k-\bar m+2},\dots, \G_{k-1}$, applying the inductive hypothesis and proceeding as in the proof of the previous theorem, we get
\begin{equation}
\begin{aligned}
|\e{ij}{k+1}|\leq \left(1-w_{min}(1-\sigma')\right)\|E^{k-l+1}\|_{\infty}
\end{aligned}
\end{equation}
with $\sigma'$ given by \eqref{eq:distltime} for distance $l-1$, 
and defining $\sigma:=\left(1-w_{min}(1-\sigma')\right)<1$ we get \eqref{eq:distltime} for  distance equal to $l$.\\
Since the sequence $\{\G_k\}$ is fully connected for sequences of length $\bar m$ we have that for every couple of nodes $i, j$ the joint distance between $j$ and $i$ in
$\G_{k-\bar m+1}, \G_{k-\bar m+2},\dots, \G_k$ is smaller or equal than $\bar m$ and we get the thesis.

\end{proof}
\end{theorem}

Lemma 1 shows that if we consider the time-independent case as the particular instance of the time-varying case where each of the graphs $\G_k$ is equal to $\G$ with diameter $\delta$, then Assumption 3 holds with $\bar m = \delta$ and the two theorems give the same inequality for the error vectors.

\section{Numerical results} 
In this section we present initial testing results for the DFIX method.
The DFIX is compared with the state-of-the-art distributed optimization algorithm from \cite{harnessing} and the method for solving systems of linear equations presented in \cite{nedic}. The test set consists of two types of problems: Simple Kriging problems and linear systems with strictly diagonally dominant co\-ef\-fi\-cient matrix. In Section 5.1 we study how the computational and com\-mu\-ni\-ca\-tion cost of DFIX is influenced by the connectivity of the underlying network and we compare DFIX with the methods from \cite{harnessing} and \cite{nedic} on a simple kriging problem. In Section 5.2 we repeat the comparison considering a randomly generated linear system. In Section 5.3 we consider the case of time-varying network.

The results demonstrate that DFIX, analogously to the classical results, outperforms the corresponding optimization method for solving the un\-con\-strained quadratic problem both in terms of computational and com\-mu\-ni\-ca\-tion costs. With respect to the method from\cite{nedic} the comparison is again favorable for DFIX, in the case of the iterative matrix with suitable properties. Clearly, the method from \cite{nedic} is designed for a wider class of problems, but its  efficiency is significantly lower than DFIX efficiency in the case of unique solution and a suitable iterative matrix.\\

For the sake of completeness we describe here both methods we compare with. We already remarked in the introduction that finding a solution of \eqref{eq:ls} is equivalent to solve the unconstrained optimization problem with quadratic objective function given by $\frac{1}{2}x^tAx - b^tx$. When applied to this optimization problem, the method from \cite{harnessing}, abbreviated as "Harnessing", can be stated as follows. 
Within one Harnessing iteration, each node computes its own solution estimate $ x_i^{k+1} $ and an ad\-di\-tion\-al vector $ s_i^{k+1}, $ which is an estimation for the average gradient, according to the following rule
\begin{eqnarray} 
x_i^{k+1} & = & \sum_{j =1}^n w_{ij} x_i^k - \eta s_i^k \label{h1}\\
s_i^{k+1}  & = & \sum_{j =1}^n w_{ij} s_i^k  + A_i (x_i^{k+1} - x_i^k) \label{h2} 
\end{eqnarray}
with $ \eta $ in (\ref{h1}) being the hand tuned step size parameter and $A_i$ denoting the $i$-th row of the matrix.

The second method \cite{nedic} we consider, abbreviated as "Projection", deals with the linear system \eqref{eq:ls} directly and is specified as follows. Before the iterative procedure starts, each agent $i$ defines the local initial vector $x^0_i$ as any solution of the equation $A_ix = b_i$ then, at every iteration, each node performs the following update:
\begin{equation*}
x^{k+1}_i = x^k_i-\frac{1}{|\mathcal{O}_i|}P_i\left(|\mathcal{O}_i|x^k_i-\sum_{j\in\mathcal{O}_i}x^k_j\right)
\end{equation*}
where $\mathcal{O}_i$ denotes the neighborhood of node $i$ in the network and $P_i$ is the projection matrix on the subspace $\ker(A_i) = \{x\in\Rn\ |\ A_ix = 0 \}$.

The DFIX method we consider here is defined using Jacobi Overrelaxation, as specified in Section 2, as underlying fixed point method. The iteration $k$ of the resulting method at each node is given by
\begin{equation}\label{eq:JORstep1}
\begin{aligned}
&\hx{ii}{k+1}=(1-\alpha)\x{ii}{k}-\frac{\alpha}{a_{ii}}\left(\sum_{j\neq i}a_{ij}\x{ij}{k}-b_i\right), \; \hx{ij}{k+1}=\x{ij}{k}\ \ \ \text{for}\  j\neq i, \\
\end{aligned}\end{equation}
and 
 \begin{equation}\label{eq:JORstep2}
\x{i}{k+1}=\sum_{j=1}^nw_{ij}\hx{j}{k+1}. \end{equation}
 In the rest of the section we refer to the method defined by equations \eqref{eq:JORstep1}, \eqref{eq:JORstep2} as DFIX - JOR.
 
\subsection{ Simple  Kriging problem}
The first problem we consider is Simple Kriging \cite{kriging}.
Let us consider  a physical process modeled as a spatial random field and assume that a network of sensors is given in the region of interest, taking measurements of the field. 
The goal is to estimate the field in any given point of the region. Assuming that the field is Gaussian and stationary, and that the expected value and covariance function are known at any point, this kind of problem can be solved by Simple Kriging method.

Denote with $\Z(s)$ the value of the random field at the point $s$, and with $\mu(s)$ its expected value, which is assumed to be known. 
Moreover, by the stationarity assumption, we have that the covariance between the value of $\Z$ at two points is given by 
$$\Cov(\Z(s_1), \Z(s_2))=\K(\|s_1-s_2\|_2)$$
for some nonnegative function $\K$.
Given $\{s_1,\dots,s_n\}\subset\R^2$ the positions in space of the $n$ sensors of the network,
let  $\{\Z(s_1),\dots,\Z(s_n)\}$ be the sampled values at those points and define the covariance matrix $A =[a_{ij}] \in\Rnn$ as 
\begin{equation*}\label{eq:krigingK}
a_{ij}=\K(\|s_i-s_j\|_2).
\end{equation*}
Now, given a point $\bar s$ where we want to estimate the field, we define the vector $b\in\Rn$ as 
\begin{equation*}\label{eq:krigingRHS}
b_{i}=\K(\|s_i-\bar s\|_2).
\end{equation*}

The predicted value of $\Z(\bar s)$ is then given by 
\begin{equation*}
\hat p(s):=\mu(\bar s)+\sum_{i=1}^n x_i(\Z(s_i)-\mu(s_i))
\end{equation*}
where $(x_1,\dots,x_n)$ is the approximate solution of the linear system
\begin{equation}\label{eq:System}
Ax=b.
\end{equation} 

Clearly, the matrix $ W $ plays an important role in the DFIX - JOR method. So let us first illustrate the influence of connectivity within the network in terms of communication traffic and computational cost for the above described kriging problem,  with covariance function given by 
\begin{equation}\label{eq:covariance}
\K(t):=\exp(-5t^2).
\end{equation}
We assume a set $\{s_1,\dots,s_{100}\}\subset[-30,30]^2$ of agents is given and for any $m\in\{2, 4, \dots, 48,50\}$ we take the $m$-regular graph with vertices $\{s_1,\dots,s_{100}\}$. That is, given the value of $m$, we define the network so that each node has degree $m$. The matrix W is defined using the Metropolis weights \cite{metropolis} which in the $m$-regular case are given by \begin{equation}\label{eq:metropolisregular}
w_{ij} = \begin{cases}
(m+1)^{-1} & \text{if}\ j=i\ \text{or}\ j\in\mathcal{O}_i\\
0 & \text{otherwise}
\end{cases}
\end{equation} For every value of the degree $m$ we apply DFIX-JOR method to solve $ Ax = b. $  

At Figure 1 and 2 we plot the number of iterations performed by the method and the total communication cost, respectively, until the stopping criterion
\begin{equation}\label{termcond}
\max_{i = 1,\ldots,n} \|Ax_i^{k} - b\| \leq 10^{-4}
\end{equation} 
is satisfied, for graphs of increasing degree.  In other words we are asking that each node solves the system with the residual tolerance of $ 10^{-4}. $ The communication cost is computed as follows. At each iteration, Step 1 does not require any communication between the agents, while in Step 2 node $i$ shares $\x{i}{k}$ with all the agents in its neighbourhood. The per-iteration traffic is thus given by $nm = 2|\E|$, where $\E$ is the set of edges of the underlying network and $m$ is the degree. \\
From Figures 1 and 2 we can see that, as the degree of the network increases, the number of iterations required to satisfy \eqref{termcond} decreases, while the total communication traffic first decreases then increases again. This behaviour can be explained as follows\footnote{Note that here we implicitly assume that there is a dedicated communication link between any pair of agents, i.e., the broadcast nature of communication is not considered. While broadcast transmissions can be considered in future studies, current comparisons are appropriate and fair and reflect practical scenarios where dedicated peer-to-peer channels are ensured, e.g., through frequency division multiple access or similar schemes.}. As the connectivity of the graph improves, the local information is distributed through the network more efficiently, and a smaller number of iterations is necessary. On the other hand, if the degree is larger, the consensus step \eqref{eq:step2} of the algorithm requires each node to share its local vector with a larger number of neighbours, yielding a higher communication traffic at each iteration. The fact that the overall communication traffic (Figure 2) is nonmonotone suggests that for large values of the degree, the decrease in the number of iterations in not enough to balance the higher per-iteration traffic.
\begin{figure}[h]
\centering
\begin{minipage}{.45\textwidth}
  \centering
  \includegraphics[width=\linewidth]{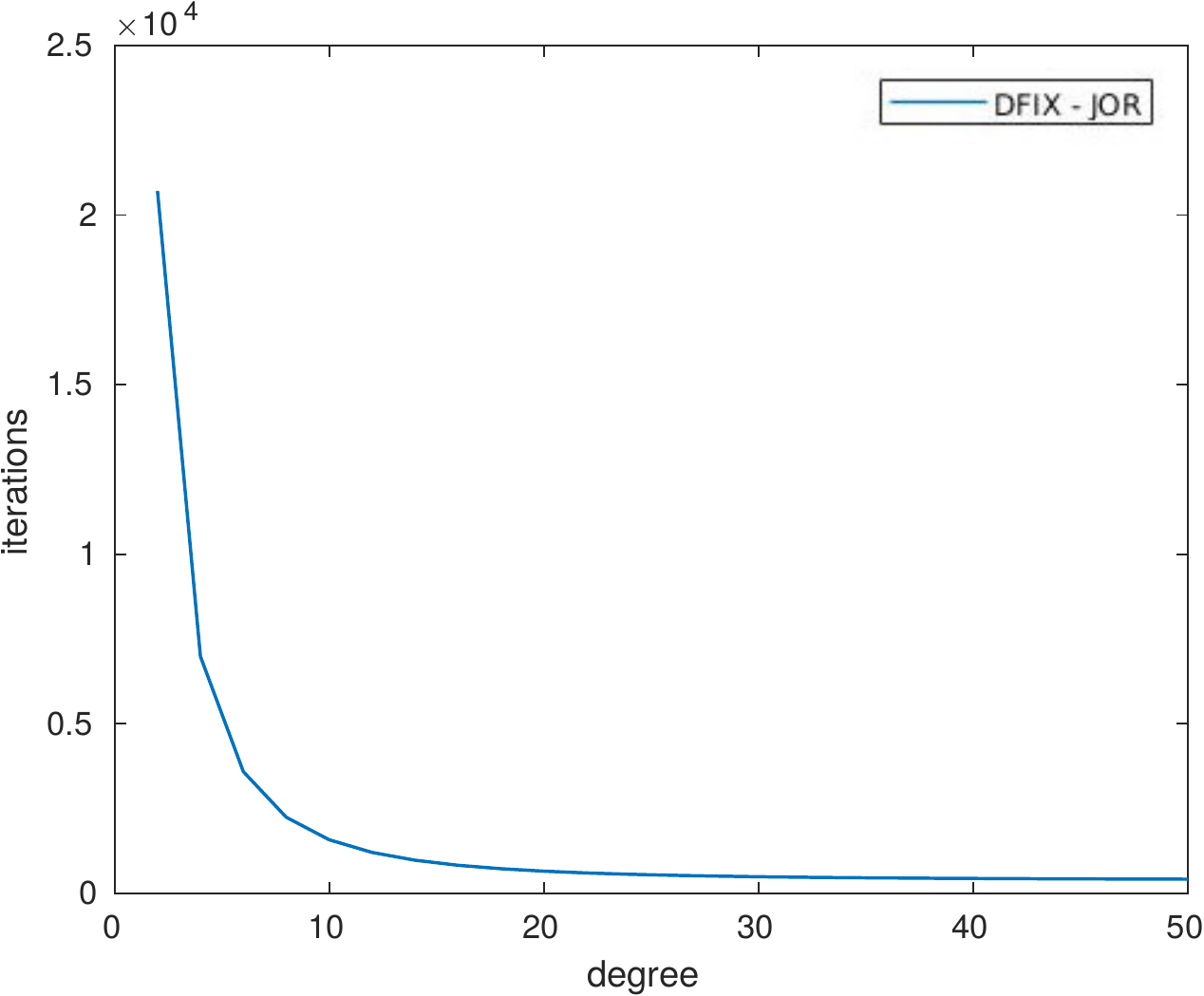}
  \caption{Number of iterations}{}
\end{minipage}%
\hfill
\begin{minipage}{.45\textwidth}
  \centering
  \includegraphics[width=\linewidth]{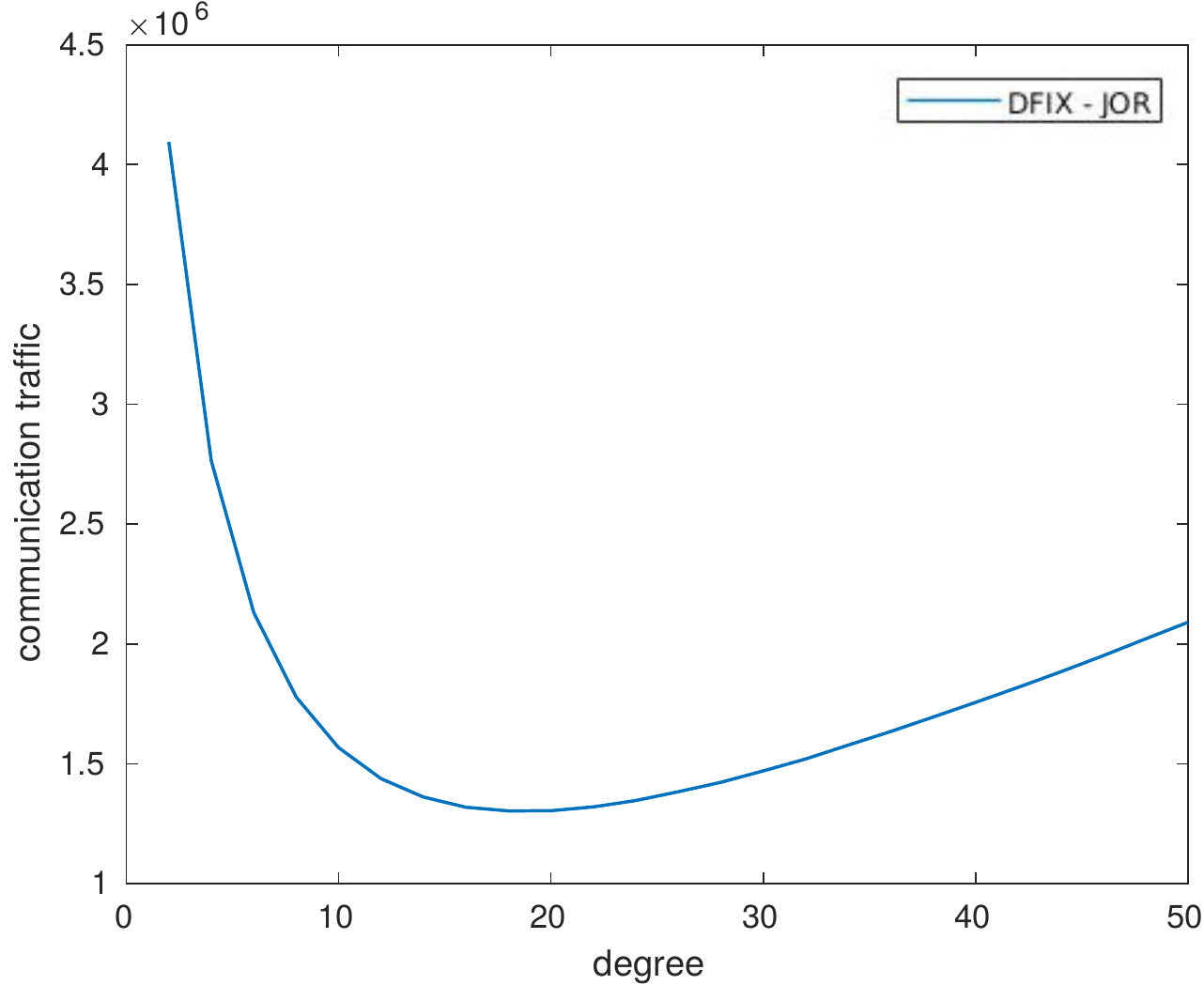}
  \caption{Communication cost}{}
\end{minipage}
\end{figure}

Let us now compare the DFIX - JOR with Harnessing \cite{harnessing} and Projection method \cite{nedic}. 
 We consider a $10\times10$ grid of nodes located at $\{s_1,\dots,s_{100}\}\subset[-3,3]^2$ and, given a communication radius $R>0$ we define the network so that nodes $i$ and $j$ are neighbours if and only if their distance is smaller than $R$. The linear system that we consider is derived by the kriging problem described at the beginning of this section. That is, we consider again $Ax=b$ with
\begin{equation} a_{ij} = \K(\|s_i-s_j\|_2),\phantom{sp} b_i = \K(\|s_i-\bar s\|_2) \label{kriging2} \end{equation}
where $\K$ is given by \eqref{eq:covariance} and $\bar s$ is a fixed random point in $[-3,3]^2$. Proceeding as in the previous test, we compute the communication traffic and com\-pu\-ta\-tion\-al cost required by the three methods to achieve the tolerance specified at \eqref{termcond}, for different values of the communication radius $R$. For each method, the overall computational cost is given by the number of iterations performed times the per-iteration cost, calculated as the number of scalar operations in one iteration. Similarly, the communication traffic is given by the number of iterations times the total number of vectors shared by the nodes during one iteration, times the length $n$ of the vector.
The matrix $W$ is defined as in \cite{metropolis}. That is, we define the off-diagonal elements as 
\begin{equation*}
w_{ij} = \begin{cases}
\frac{1}{1+\max\{m_i, m_j\}} & \text{if}\ j\in\mathcal O_i\\
0 & \text{otherwise}
\end{cases}
\end{equation*}
where $m_i$ denotes the degree of node $i$, and the diagonal elements as 
$$w_{ii} = 1-\sum_{j\neq i} w_{ij}$$ 
so that the resulting matrix $W$ is stochastic.
The stopping criterion is the same as in the previous test, i.e., each node solves the problem with the tolerance of $ 10 ^{-4}. $ The initial point at each node is the same for the three methods and is defined as follows:
$$x^0_{ii} = \frac{b_i}{a_{ii}}, \phantom{sp} x^0_{ij} = 0\ \text{for every}\ j\neq i. $$
Moreover, the relaxation parameter $\alpha$ in \eqref{eq:JORstep1} is chosen as $\frac{2}{\|D^{-1}A\|_\infty}$ where $D=\diag(a_{11},\ldots, a_{nn}), $
while for Harnessing method we take in \eqref{h1}
$\eta = \frac{1}{3L}$ where $L =\max_{i=1:n} 2\|A_i\|_2^2$.\\

In Figures 3 and 4 we plot the obtained results. As we can see, in this framework, DFIX method is more efficient than the two methods we compare with, both in terms of computational cost and in terms of communication traffic.
$\ $\\
\begin{figure}[h]
\begin{minipage}{.45\textwidth}
  \centering
  \includegraphics[width=\linewidth]{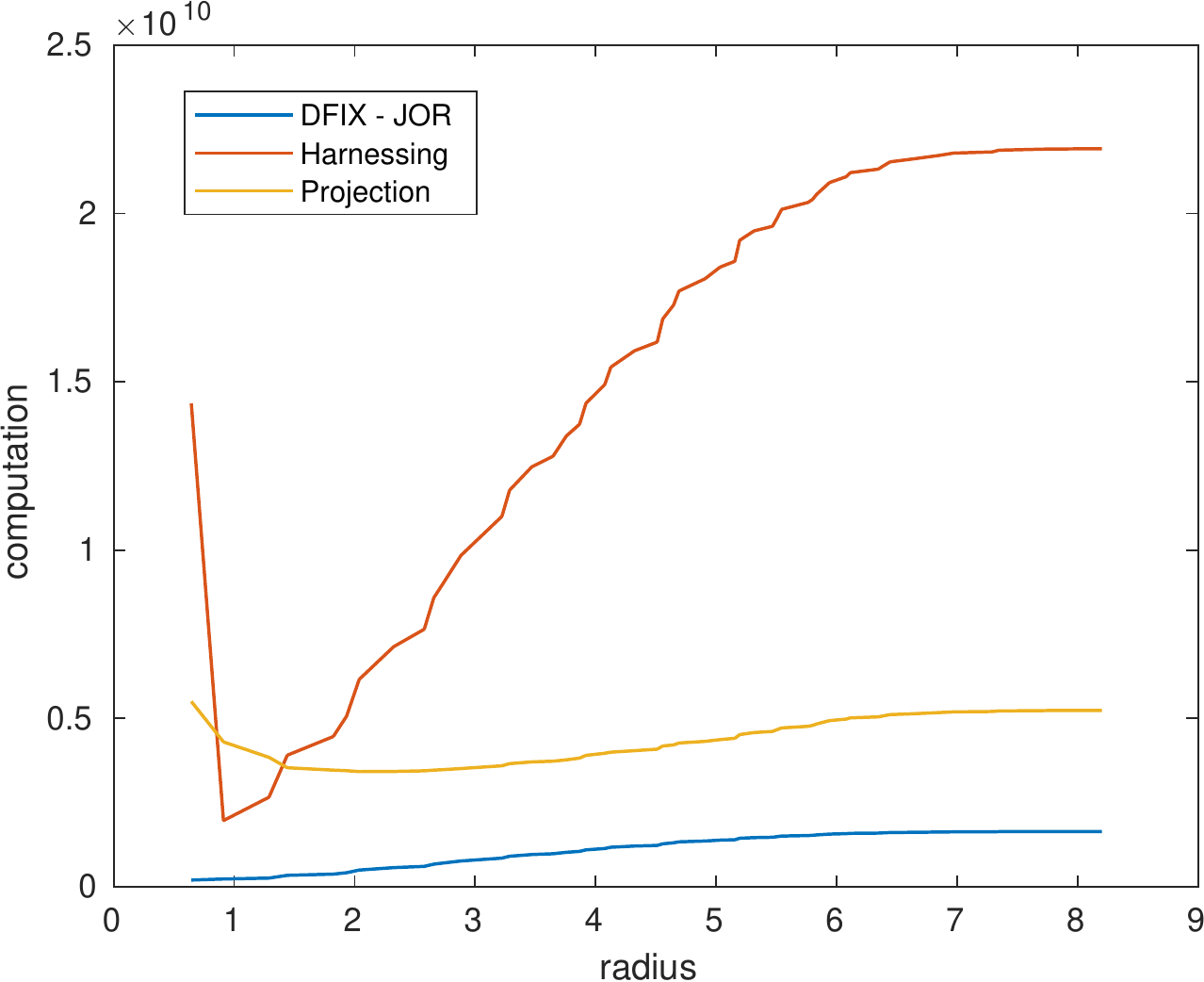}
  \caption{Simple kriging problem (\ref{kriging2}), computational cost}{}
\end{minipage}
\hfill
\begin{minipage}{.45\textwidth}
  \centering
  \includegraphics[width=\linewidth]{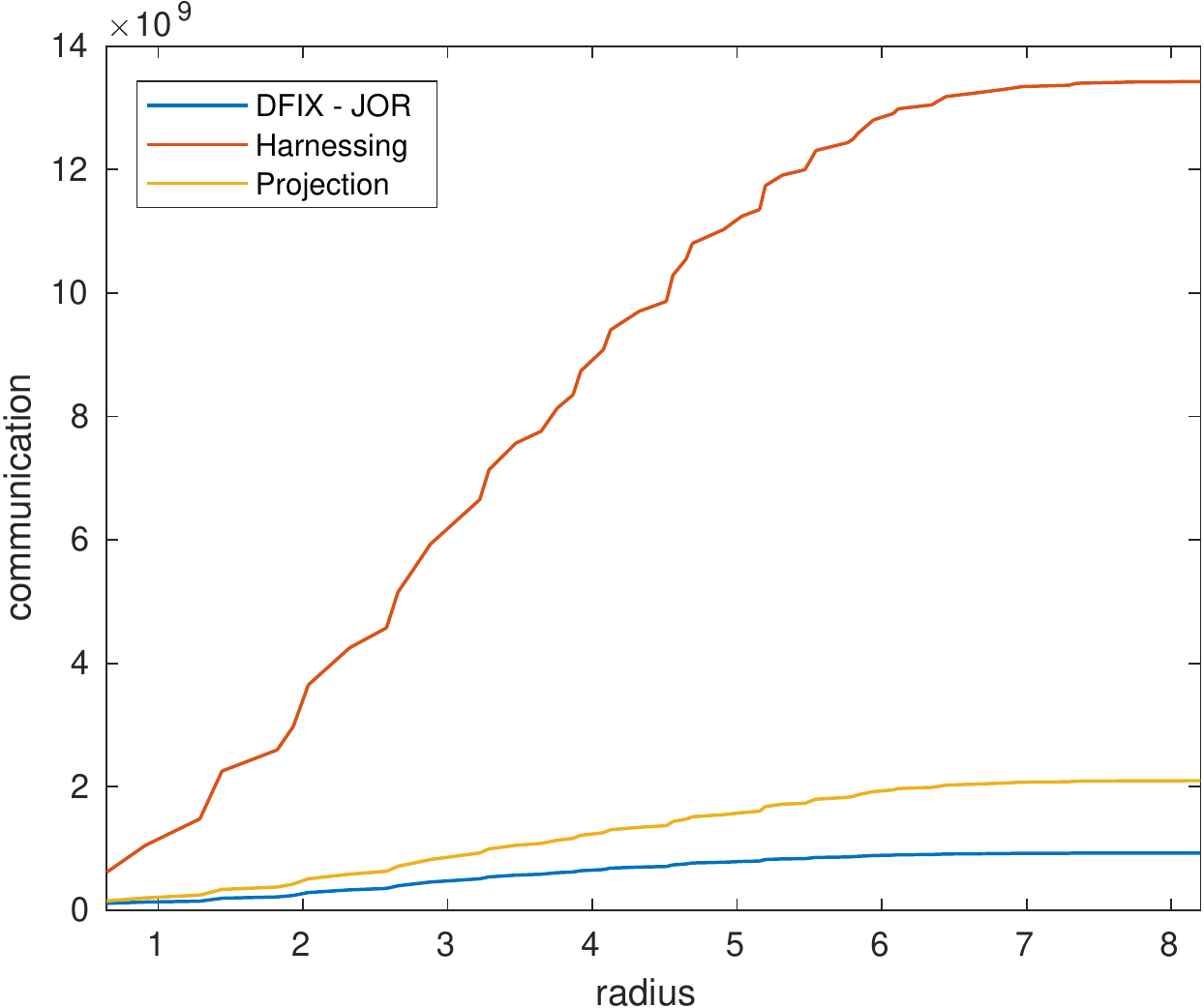}
  \caption{Simple kriging problem (\ref{kriging2}), communication traffic}{}
\end{minipage}
\end{figure}

\subsection{Strictly diagonally dominant systems}
Let us now  consider a linear system $Ax=b$ of order $n = 100$, where $ A $ and $ b $ are generated as follows. For every index $i$ we take $b_i$ randomly generated with uniform distribution in $(0,1)$, and $A$ is a symmetric diagonally dominant random matrix obtained as follows: take $\hat{a}_{ij}\in(0,1)$ with uniform distribution and  then set $\tilde{A}=\frac{1}{2}(\hat{A}+\hat{A}^T)$ and finally $A = \hat{A}+(n-1)I$, where we denote with $I$ the identity matrix of order $n$. 
 As the underlying network we consider an $m$-regular graph with $n$ nodes. For every fixed value of the degree $m$ we generate, as just described, 10 random linear systems, solve all of them using the three methods and compute the average number of iterations necessary to arrive at termination. For each method, the total amount of computation and communication are then obtained multiplying the average number of iterations and the per-iteration computational cost and communication traffic, respectively.
The matrix $W$ is defined as in \eqref{eq:metropolisregular}, the step sizes $\alpha$ and $\eta$, the initial guess at each node and the termination condition are as in the previous test.
In Figures 5 and 6 we plot the results  for $m$ in in $\{2, 4, \dots, 48, 50\}$.
Similarly to the previous test, we have that DFIX outperforms both Harnessing an Projection method in terms of computation and communication. From Figure 6 we can notice that the communication required by the two methods for distributed linear systems, DFIX and Projection, is similar and that the difference with the communication required by Harnessing method increases as the degree of the graph increases. Regarding the computational cost (Figure 5), we have that while DFIX is cheaper than the other two methods, Projection method seems to be more influenced by the connectivity of the network and it is more efficient than Harnessing only for large vaues of the degree.\\

\begin{figure}[h]
\begin{minipage}{.45\textwidth}
  \centering
  \includegraphics[width=\linewidth]{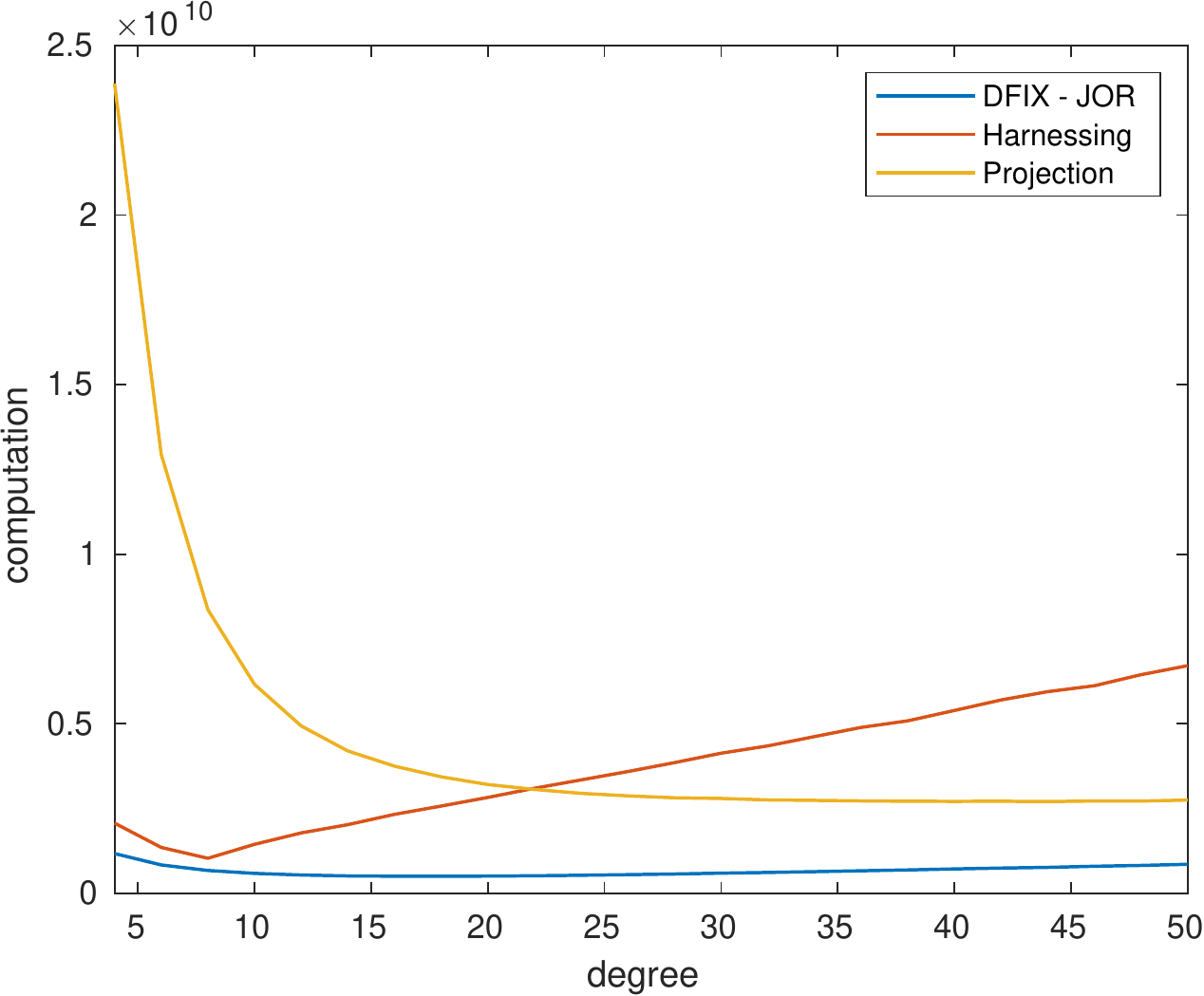}
  \caption{$m$-regular graph, computational cost}{}
\end{minipage}
\hfill
\begin{minipage}{.45\textwidth}
  \centering
  \includegraphics[width=\linewidth]{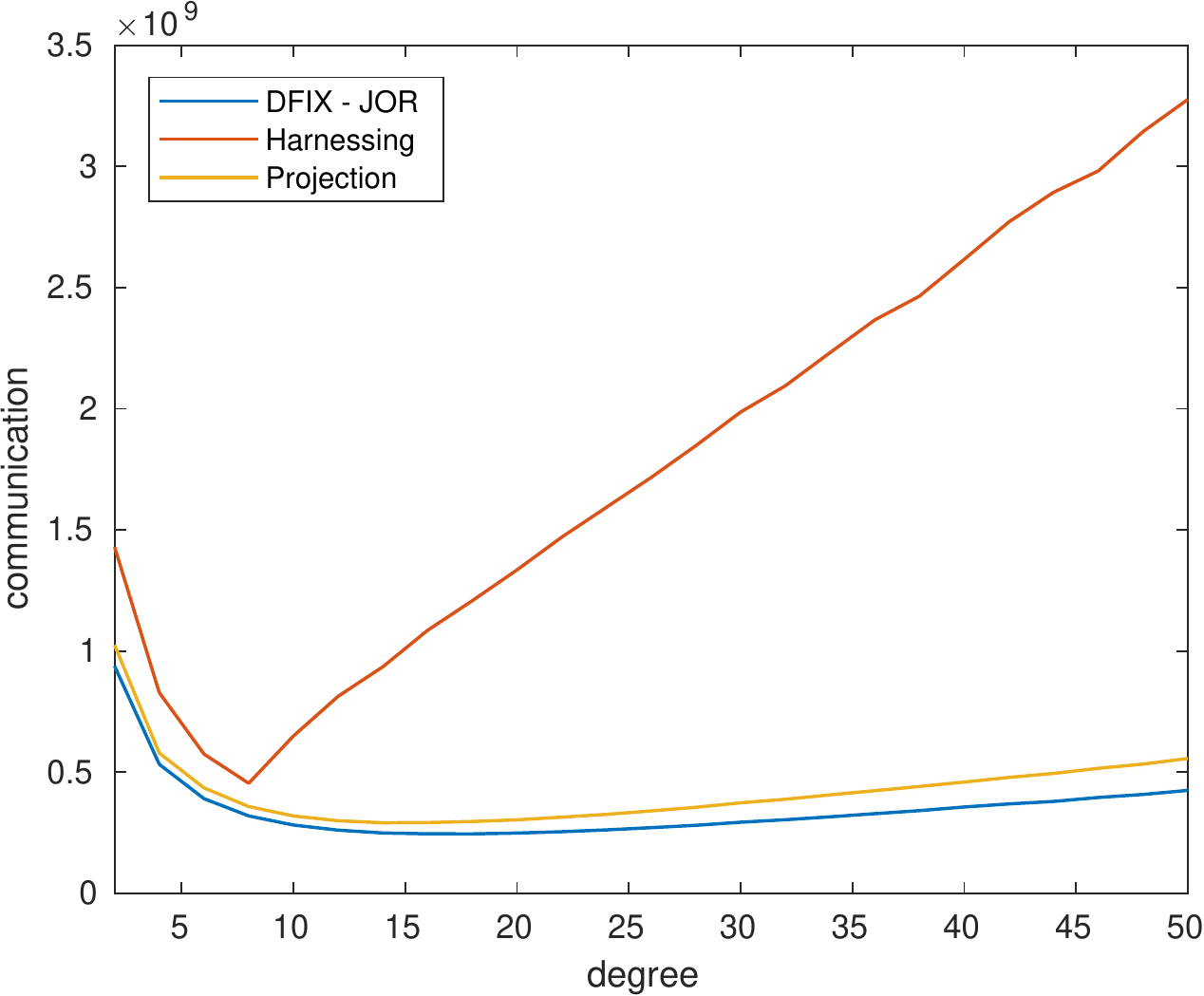}
  \caption{$m$-regular graph, communication traffic}{}
\end{minipage}
\end{figure}

\subsection{Time-varying Network}

We now compare the performance of the three methods in the time-varying case described in Section 4. 

The sequence $\{\G_k\}$ is generated as follows. We first fix a strongly connected graph $\G = (\V, \E)$ and a scalar $\gamma\in(0,1]$. Then, at every iteration $k$ we randomly generate $\E_k$ by uniformly sampling $\gamma |E|$ edges from $\E$ and we set $\G_k = (\V, \E_k)$. This construction can be interpreted as having a fixed underlying graph $\G$ that represents the available communication links among the nodes, and employing at each iteration only a fraction $\gamma$ of the links. In particular, $\gamma = 1$ corresponds to the case when $\G_k = \G$ for every $k$. As remarked in Section 4, this is equivalent to the time-independent case.

\begin{figure}[h]
\begin{minipage}{.45\textwidth}
  \centering
  \includegraphics[width=\linewidth]{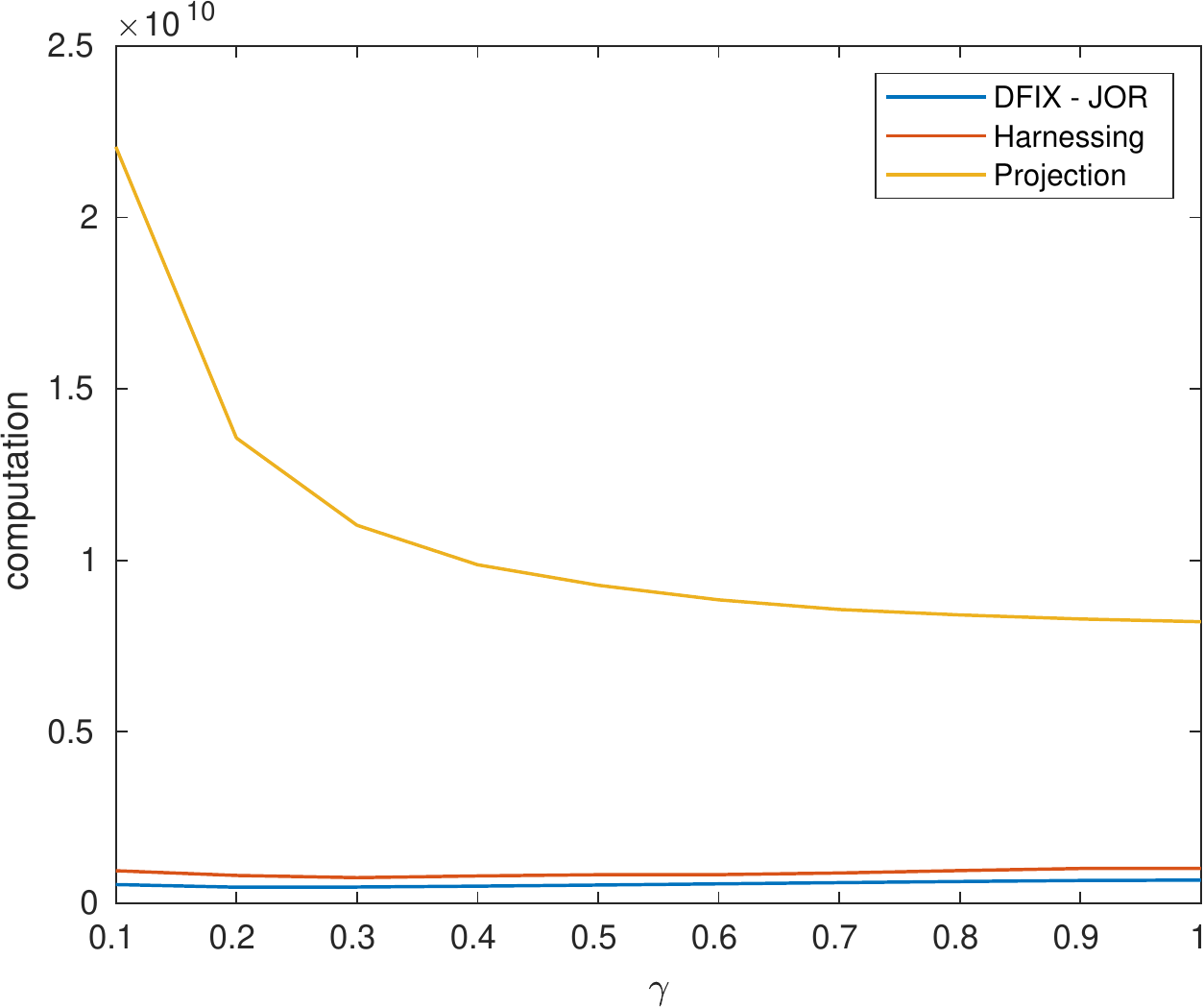}
  \caption{computational cost}{}
  \label{fig:time_comp}
\end{minipage}
\hfill
\begin{minipage}{.45\textwidth}
  \centering
  \includegraphics[width=\linewidth]{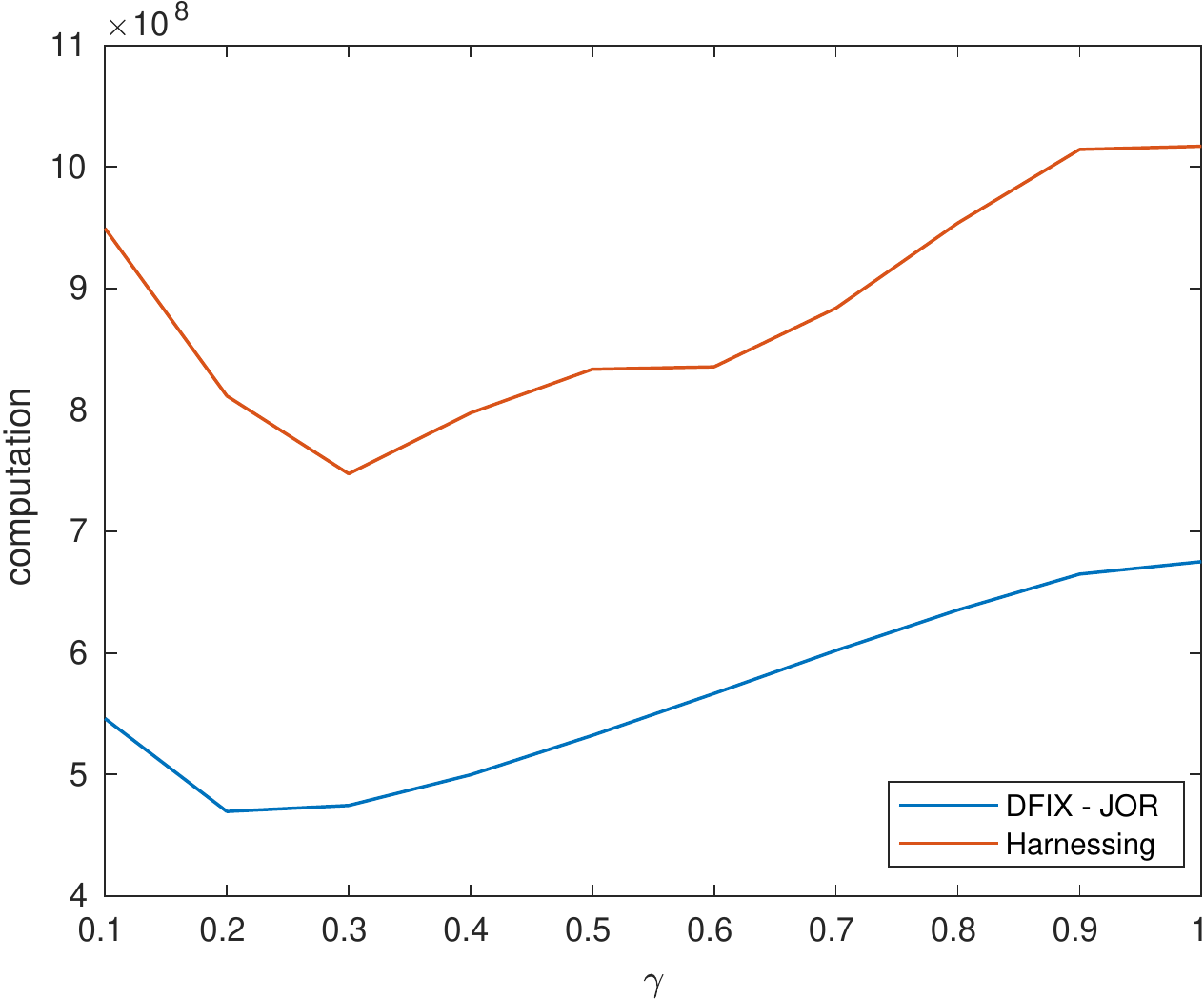}
  \caption{computational cost}{}
  \label{fig:time_compDH}
\end{minipage}\\

 \centering
\begin{minipage}{.45\textwidth}
  \includegraphics[width=\linewidth]{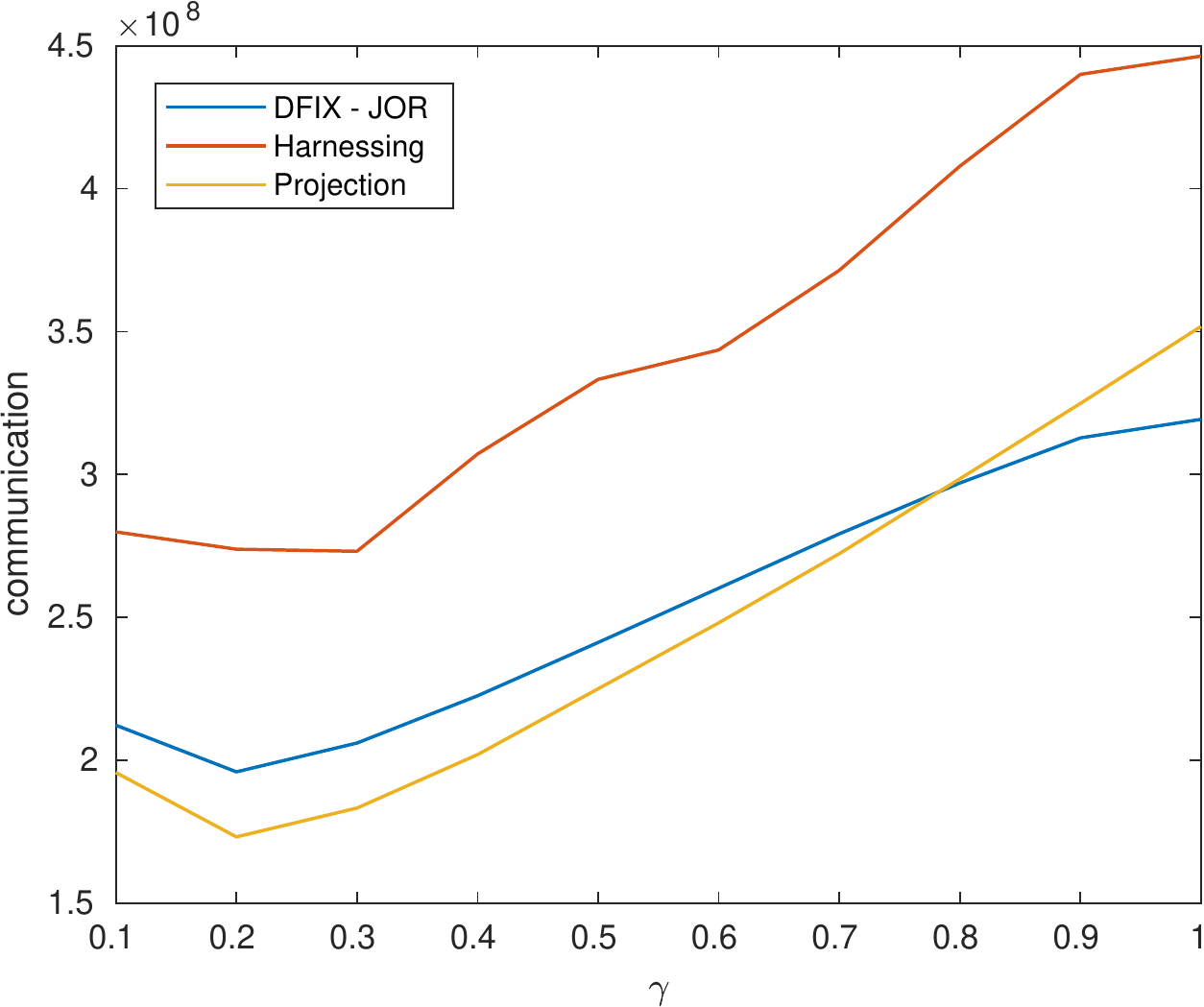}
  \caption{communication traffic}{}
  \label{fig:time_comm}
\end{minipage}
\end{figure}

The test we present here is carried on comparing the communication and computational cost required by the three methods to solve a given linear system using the same sequence of networks $\{\G_k\}$. In practice we generated the linear system as in Section 5.2 and we chose $\G$ as the undirected $m$-regular graph with $n=100$ vertices and degree $m=8$. We repeated the same test for $\gamma$ in $\{0.1, 0.2, \dots, 1\}$. For every $k$ the consensus matrix $W^k$ associated with $\G_k$ is defined as in \eqref{eq:metropolisregular}, the terminantion condition and all the prameters of the methods are chosen as in the previous sections. In Figures \ref{fig:time_comp}, \ref{fig:time_compDH} and \ref{fig:time_comm} we plot the results (Note that Figure 8 repeats the results of Figure 7, focusing only on the comparison Harnessing versus DFIX-JOR). The computational cost and the communication traffic are calculated as described in Section 4.2.  

We can see that, in the considered framework, DFIX outperforms Har\-ness\-ing method both in terms of computation and communication. Com\-pa\-ring with Projection, we have that, for every value of the parameter $\gamma$, the computational cost of DFIX is significantly lower, but it requires a smaller amount of communication only for large values of $\gamma$ (that is, when each graph $\G_k$ is equal or close to $\G$). Moreover we can see that for all the methods there is an optimal value of $\gamma<1$, that minimizes the communication traffic, suggesting that using the whole graph $\G$ at every iteration (that is, setting $\gamma=1$) is unefficient. A similar phenomena happens for Harnessing and DFIX also for the computational cost (Figure 8), while we can see in Figure 7 that Projection method is most efficient when all the available communication links are used at each iterations. For $\gamma<1$ the networks $\G_k$ are in general not connected, but the joint connectivity of the overall sequence is enough to ensure the convergence of the methods.

\section{Conclusions}
We proposed a class of novel, iterative, distributed methods for the solution of linear systems of equations, derived upon classical fixed point methods. We proved global convergence in the case when the communication network is strongly connected and we showed that the convergence rate depends on the diameter of the network and on the norm of the underlying iterative matrix. In particular we have that if the graph is strongly connected, the obtained result is analogous to the classical, centralized, case. We extended the presented method to the time-varying case and we proved an analogous convergence result, assuming the networks satisfy suitable joint con\-nec\-tivity assumptions, comparable with those required by different methods in literature.

Our algorithm was compared with the relevant methods presented in \cite{harnessing} and \cite{nedic}. The numerical results showed good performance of DFIX compared with the mentioned methods. In particular, in the vast majority of the considered tests, DFIX outperformed the two methods in terms of both computational cost and communication traffic.

\section*{Acknowledgements}
This work is supported by the BIGMATH project which has received funding by the European Union’s Horizon 2020 research and innovation programme under the Marie Sk\l{}odowska-Curie Grant Agreement no. 812912.
The work of Jakoveti\'c, Kreji\'c and Krklec Jerinki\'c is  partially supported by Serbian Ministry of Education, Science and Technological Development, grant no. 174030.


\begin{thebibliography}{99}
\bibitem{Berman}A. Berman, R. J. Plemmons, Nonnegative Matrices in the Mathematical Sciences, SIAM, 1994

\bibitem{bertsekas} D.P. Bertsekas, J. N. Tsitsiklis, Parallel and Distributed Computation: Numerical Methods, Athena Scientific 2015. 

\bibitem{kriging} N. A. C. Cressie, Statistics for Spatial Data. New York: Wiley, 1993

\bibitem{consensus} M. H. deGroot, Reaching a Consensus,  J.  Amer.  Statist.  Assoc.,  vol. 69, no. 345, pp. 118-121, 1974.

\bibitem{fm} A. Frommer, G. Mayer, Convergence of Relaxed Parallel Multisplitting Methods, Linear Algebra Appl., 119 (1989), pp. 141-152 

\bibitem{fs1} A. Frommer, D.B. Szyld,  H-splittings and Two-Stage Iterative Methods, Numer. Math., 63 (1992), pp. 345-356
 
\bibitem{fs2} A. Frommer, D.B. Szyld, On Asynchronous Iterations, Journal of Computational and Applied Mathematics. 123 (2000), 201-216. 

\bibitem{greenbaum} A. Greenbaum, Iterative Methods for Solving Linear Systems, SIAM, (1997).

\bibitem{graphs} J. M. Hendrickx, R. M. Jungers, A. Olshevsky, G.  Vankeerberghen, Graph Diameter, Eigenvalues, and Minimum-Time Consensus, Automatica 50 (2014), 635-640 .
\bibitem{krige} D. G. Krige, A Statistical Approach to Some Basic Mine Valuation Problems on the Witwatersrand: J. Chem. Metal. Min. Soc. South Africa, v. 52 (1951), p. 119–139.

\bibitem{harnessing} N. Li, G. Qu, Harnessing Smothness to Accelerate Distributed Optimization, IEEE Transactions  Control of Network Systems 5 (3), (2017), 1245-1260

\bibitem{nedic} Ji Liu, A. S. Morse, A. Nedi\'c, T. Basar, Exponential convergence of a distributed algorithm for solving linear algebraic equations, Automatica 83 (2017), 37-46.  

\bibitem{dls2} J. Liu , S. Mou, A. S. Morse, Asynchronous Distributed Algorithms for Solving Linear Algebraic Equations, IEEE Transactions on Automatic Control,  Vol. 63, No. 2, (2018), 372-385. 

\bibitem{matheron} G. Matheron, Traité de geostatistique appliquée, vol. II, Le krigeage: Memoires du Bureau de Recherches Géologiques et Miniéres, no. 24, Editions Bureau de Recherche Géologiques et Miniéres, Paris, 1963.

\bibitem{dls1} S. Mou, Z. Lin, L. Wang, D. Fullmer, A.S. Morse, A Distributed Algorithm for Efficiently Solving Linear Equations and Its Applications (Special Issue JCW), System \& Control Letters, 91 (2016), 21-27. 

\bibitem{diging} A. Nedic, A. Olshevsky, W. Shi, Achieving geometric convergence
for distributed optimization over time-varying graphs, SIAM J. Optim., 27 (2016).

\bibitem{saad} Y. Saad, Iterative Methods for Sparse Linear Systems, SIAM, 2003

\bibitem{extra} W. Shi, Q. Ling, G. Wu, W. Yin, EXTRA: An Exact First-Order Algorithm for Decentralized Consensus Optimization, SIAM J. Optim., 25 (2015), 944–966

\bibitem{touri} B. Touri and A. Nedic, On Backward Product of Stochastic Matrices, Automatica 48 (2012), 1477-1488.

\bibitem{survey} P. Wang, S. Mou, J. Lian, W. Ren, Solving a System of Linear Equations: From Centralized to Distributed Algorithms, Annual Reviews in Control, (to appear)


\bibitem{dls3} X. Wang, J. Zhou, S. Mou, M. J. Corless, A Distributed Algorithm for Least Square Solutions of Linear Equations, https://arxiv.org/pdf/1709.10157.pdf

\bibitem{dls4} Y. Xiao, J. Hu, Distributed Solutions of Convex Feasibility Problems with Sparsely Coupled Constraints, 2017 IEEE 56th Annual Conference on Decision and Control (CDC), 3386-3392

\bibitem{metropolis} L. Xiao, S. Boyd, and S. Lall, Distributed average consensus with time-varying Metropolis weights, 2006























\end{thebibliography}
\end{document}